\documentclass{article}
\title{There are no socialist primes less than $10^{9}$}
\author{
Tim Trudgian\footnote{Supported by Australian Research Council DECRA Grant DE120100173.}\\
Mathematical Sciences Institute\\ The Australian National University,
 ACT 0200, Australia\\ timothy.trudgian@anu.edu.au
}
\usepackage{url}
\usepackage{amsthm}
\usepackage{amsmath}

\usepackage{amssymb}
\usepackage{booktabs}

\begin{document}

\maketitle
\begin{abstract}
\noindent
There are no primes $p$ with $5<p<10^{9}$ for which $2!, 3!, \ldots, (p-1)!$ are all distinct modulo $p$; it is conjectured that there are no such primes.
\end{abstract}

\section{The problem}
Erd\H{o}s asked whether there exist any primes $p>5$ for which the numbers $2!, 3!, \ldots, (p-1)!$ are all distinct modulo $p$. Were these $p-2$ factorials all distinct then the $p-1$ non-zero residue classes modulo $p$ contain at most one of them. Motivated by this redistribution of resources amongst classes I shall call such a prime $p$ a \textit{socialist prime}.

Rokowska and Schinzel \cite{RokS}\footnote{This problem also appears as \textbf{F11} in Richard Guy's insuperable book \cite{GuyBook}.}
 proved that $p$ is a socialist prime only if $p \equiv 5 \pmod 8$, and 
\begin{equation}\label{con}
\left(\frac{5}{p}\right) = -1, \quad \left( \frac{-23}{p}\right) = 1.
\end{equation}
Moreover, if a socialist prime exists then none of the numbers $2!, 3!, \ldots, (p-1)!$ is congruent to $-((p-1)/2)!$. 
The proof given by Rokowska and Schinzel is fairly straightforward. 

One may dismiss primes of the form $p \equiv 3 \pmod 4$, since such primes have the property \cite[Thm 114]{HW} that $((p-1)/2)! \equiv \pm 1 \pmod p$. By Wilson's theorem, $(p-1)! \equiv -1 \pmod p$ and  $(p-2)! \equiv (p-1)! (p-1)^{-1} \equiv +1 \pmod p$, conditions which, when taken together, prohibit $p$ from being a socialist prime. Henceforth consider $p \equiv 1 \pmod 4,$ in which case
\begin{equation}\label{gg}
\left\{ \left( \frac{p-1}{2}\right)! \right\}^{2} \equiv -1 \pmod p.
\end{equation}

If $2!, 3!, \ldots, (p-1)!$ are all distinct modulo $p$ then they must be permutations of the numbers $1,2,\ldots, p-1$ with the exception of some $r$, with $1 \leq r \leq p-1$, whence
$$\prod_{n=2}^{p-1} n! \equiv \frac{(p-1)!}{r} \pmod p,$$
so that 
$$1 \equiv r \prod_{n=1}^{p-2} n! \equiv r((p-1)/2)! \prod_{1 \leq k < \frac{p-1}{2}} k! (p-k-1)! \pmod p.$$
Applying (\ref{gg}) and Wilson's theorem gives\begin{equation*}\label{g2}
r \prod_{1 \leq k < \frac{p-1}{2}} (-1)^{k+1}  \equiv - \left( \frac{p-1}{2}\right) ! \pmod p,
\end{equation*}
so that $r \equiv \pm ((p-1)/2)! \pmod p$. One may dismiss the positive root, since $r$ is not congruent to any $j!$ for $1\leq j \leq p-1$. Hence 
\begin{equation*}
\prod_{1 \leq k < \frac{p-1}{2}} (-1)^{k+1} \equiv 1 \pmod p.
\end{equation*}
Equating powers of $(-1)$ gives
$$ \sum_{1\leq k < \frac{p-1}{2}} (k+1) = \frac{(p-3)(p+3)}{8} \equiv 0 \pmod 2,$$
whence, since $p\equiv 1 \pmod 4$, one may conclude that $p \equiv 5 \pmod 8.$

The conditions in (\ref{con}) are a little more subtle. Consider a polynomial $F(x) = x^{n} + a_{1} x^{n-1} + \ldots + a_{0}$ with integral coefficients and discriminant $D$. A theorem by Stickelberger (see, e.g.\ \cite[p.\ 249]{Dickson}) gives $\left( \frac{D}{p}\right) = (-1)^{n - \nu}$, where $\nu$ is the number of factors of $F(x)$ that are irreducible modulo $p$.
Consider the two congruences 
$$x(x+1) -1 \equiv 0 \pmod p, \quad x(x+1)(x+2) -1 \equiv 0 \pmod p,$$
the polynomials in which have discriminants 5 and $-23$. 
For the former, if $\left(\frac{5}{p}\right) = 1$, then, by Stickelberger's theorem, there are two irreducible factors, whence the congruence factors and has a solution. Therefore $(x+1)! \equiv (x-1)! \pmod p$ and $p$ is not a socialist prime.
Likewise for the latter: if $\left(\frac{-23}{p}\right) = -1$ then there are two irreducible factors, whence $(x+2)! \equiv (x-1)! \pmod p$.

One cannot continue down this path directly. Consider $x(x+1)(x+2)(x+3) -1 \equiv 0 \pmod p$ which has a solution if and only if $y(y+2) -1 \equiv 0 \pmod p$ has a solution, where $y = x(x+3)$. Hence $(y+1)^{2} \equiv 2 \pmod p$, which implies 2 is a quadratic residue modulo $p$ --- a contradiction since $p \equiv 5 \pmod 8$. 

Instead one can consider the congruence

$$x(x+1)(x+2)(x+3)(x+4)(x+5) -1 \equiv 0 \pmod p,$$
which is soluble precisely when $y(y+4)(y+6) -1\equiv 0 \pmod p$ is soluble, where $y = x(x+5)$. The cubic congruence in $y$ has discriminant $1957$, whence, by Stickelberger's theorem, if $\left(\frac{1957}{p}\right) = -1$ then $y(y+4)(y+6)$ has a linear factor. To deduce that  $(x+5)! \equiv (x-1)! \pmod p$ we need to know that
 $y \equiv x(x+5) \pmod p$ is soluble, that is, we need to know that $4y + 25$ is a quadratic residue modulo $p$. We can therefore add a condition to (\ref{con}), namely, a necessary condition that $p$ be a socialist prime is 
\begin{equation}\label{extrac}
\begin{split}
\left(\frac{1957}{p}\right) &= 1, \quad \textrm{or}\\
\left(\frac{1957}{p}\right) &= -1 \quad \& \quad \left(\frac{4y + 25}{p}\right) = -1,\\
&\qquad\qquad\qquad\textrm{for all $y$ satisfying}\quad y(y+4)(y+6) -1 \equiv 0 \pmod p
\end{split}
\end{equation}

\section{Computation and conclusion}
Rokowska and Schinzel showed that the only primes $5<p<1000$ satisfying $p \equiv 5 \pmod 8$ and (\ref{con}) were
$$13, 173, 197, 277, 317, 397, 653, 853, 877, 997.$$
Using Jacobi's \textit{Canon arithmeticus} they showed that for each prime there existed $1<k<j\leq p-1$ for which $k! \equiv j! \pmod p$.

I am grateful to Dr David Harvey who extended this to show that there are no socialist primes less than $10^{6}$. This computation took 45 minutes on a 1.7 GHz Intel Core i7 machine. Professor Tom\'{a}s Oliveira e Silva extended this to $p< 10^{9}$, a calculation which took 3 days.


The following example shows the utility of adding the condition (\ref{extrac}). Using the conditions $p\equiv 5 \pmod 8$ and (\ref{con}), it is easy to check that there are at most $4908$ socialist primes up to $10^{6}$. These need to be checked to see whether there are values of $k$ and $j$ for which $k! \equiv j! \pmod p$. Including the condition (\ref{extrac}) means that there are at most $3662$ socialist primes up to $10^{6}$ that need to be checked.

To extend the range of computation beyond $10^{9}$ it would be desirable to add another condition arising from a suitable congruence. The congruence leading to (\ref{extrac}) was of degree 6; no other suitable congruence was found for degrees 8 and~9.

In \cite{BanksLucaSS} the authors consider $F(p)$ defined to be the number of distinct residue classes modulo $p$ that are not contained in the sequence $1!, 2!, 3!, \ldots$. They show that $\limsup_{p\rightarrow\infty} F(p) = \infty$; for the problem involving socialist primes one wishes to show that $F(p)= 2$ never occurs. It would therefore be of interest to study small values of $F(p)$.

Finally, one may examine the problem na\"{\i}vely as follows. Ignore the conditions $p \equiv 5 \pmod 8$ and (\ref{con}) --- including these only reduces the likelihood of there being socialist primes. For $2\leq k \neq j \leq p-2$ we want $p\nmid j! - k!$. There are $\binom{p-3}{2} = (p-3)(p-4)/2$ admissible values of $(k, j)$. Assuming, speciously,  that the probability that $p$ does not divide $N$ `random' integers is $(1-1/p)^{N}$ one concludes that the probability of finding a socialist prime is
\begin{equation*}
\left( 1 - \frac{1}{p}\right)^{\frac{(p-3)(p-4)}{2}} \rightarrow e^{\frac{(7-p)}{2}},
\end{equation*}
for large $p$.

Given this estimate, and the computational data, it seems reasonable to conjecture that there are no socialist primes.

\section*{Acknowledgements}
I am grateful to David Harvey and Tom\'{a}s Oliveira e Silva for their computations, and to Victor Scharaschkin and Igor Shparlinski for their comments and suggestions.

\bibliographystyle{plain}
\bibliography{themastercanada}

\end{document}